\documentclass{article}

 \UseRawInputEncoding 

\newtheorem{theorem}{Theorem}[section]
\newtheorem{lemma}[theorem]{Lemma}
\newtheorem{corollary}[theorem]{Corollary}
\newtheorem{definition}[theorem]{Definition}
\newtheorem{problem}[theorem]{Problem}

\newcommand\prob{{\mbox Prob}}
\newcommand\php{{\mbox{PHP}}}
\newcommand\pphp{{\neg{{\cal P}{\cal H}{\cal P}}}}
\newcommand\fa{{\cal F}}
\newcommand\ea{{\cal E}}
\newcommand\ff{{\mathbf F}}
\newcommand\fp{{{\ff}_p}}
\newcommand\rx{{\fp[\overline x]}}
\newcommand\rr{{\fp[\overline x, \overline r]}}
\newcommand\re{{\cal R}({\cal E})}
\newcommand\rd{{\fp^{\le d}[\overline x]}}

\newcommand\phpd{{\fp^{\le d}[Var({\pphp}_n)]}}

\newcommand\one{{\Omega_{n,\epsilon}}}
\newcommand\fn{{\cal F}_n}

\newcommand\qed{\begin{flushright} {\bf q.e.d.} \end{flushright} }
\newcommand\prf{\noindent {\bf Proof :}}  
\newcommand\bits{\{0,1\}}

\begin{document}

\title{Extended Nullstellensatz proof systems}

\author{Jan Kraj\'{\i}\v{c}ek}

\date{Faculty of Mathematics and Physics\\
Charles University\thanks{
Sokolovsk\' a 83, Prague, 186 75,
The Czech Republic, {\tt krajicek@karlin.mff.cuni.cz}}}

\maketitle

\begin{abstract}
	
For a finite set $\fa$ of polynomials from $\rx$ ($p$ is a fixed prime) containing all polynomials $x^2 - x$, a Nullstellensatz proof of the unsolvability of the system
$$
	f = 0\ ,\ \mbox{ all } f \in {\cal F}
$$
in $\fp$ is an $\rx$-linear combination $\sum_{f \in {\cal F}} \ h_f \cdot f$
that equals to $1$ in $\rx$. The measure of complexity of such a proof is its degree:
$\max_f deg(h_f f)$.

We study the problem to establish degree lower bounds for some {\em extended}
NS proof systems: these systems prove the unsolvability of $\fa$ (in $\fp$)
by proving the unsolvability
of a bigger set $\fa\cup \ea$, where the set 
$\ea \subseteq \rr$ contains all polynomials $r^p - r$ and
satisfies the following soundness condition:
\begin{itemize}
	\item Any $0,1$-assignment $\overline a$ to variables $\overline x$
	can be appended by an $\fp$-assignment $\overline b$ to variables  
	$\overline r$ such that for all $g \in 	\ea$ it holds that $g(\overline a, \overline b) = 0$.
\end{itemize}
We define a notion of pseudo-solutions of $\fa$ and prove that the existence of 
pseudo-solutions with suitable parameters
implies lower bounds for two extended NS proof systems ENS and UENS defined in \cite{BIKPRS}. Further we give a combinatorial example of $\fa$ and candidate pseudo-solutions based on the pigeonhole principle.
\end{abstract}

\section{Introduction}

Let $\fp$ be the $p$-element field where $p$ is a fixed prime,  and let
$\fp[\overline x]$ be the ring
of polynomials with variables $\overline x = x_1, \dots, x_N$, for some $N \geq 1$.
Given a finite set ${\cal F} \subseteq \rx$ 
containing all polynomials $x_i^2 - x_i$, $1 \le i \le N$,
we consider the question whether the 
system of polynomial equations
\begin{equation} \label{8.1.23a}
	f = 0\ ,\ \mbox{ all } f \in {\cal F}
\end{equation}
has a solution
(we shall simply speak about the solvability of $\fa$ rather than of (\ref{8.1.23a})). 
Due to the presence of polynomials $x_i^2 - x_i$ such a solution must be
Boolean, i.e. from $\{0,1\}$. 
The restriction to Boolean values is due to our primary interest in proof complexity.

In computational and proof complexity we need to represent polynomials as input to algorithms.
We shall use dense representation: a list of coefficients from $\fp$ 
for all monomials up to the degree of the respective polynomial. Hence a degree $d$ polynomial over $N$ variables is represented by a string of $O(N^d)$ bits. 
Multiplying or adding polynomials or checking whether two 
polynomials are identical can be then done
by a polynomial time (shortly p-time) algorithm. Note that the bit size of $\fa$ is $O(N^d |\fa|)$ if $d$ bounds the degree of all polynomials in $\fa$.

The problem of solvability of $\fa$ in $\fp$
belongs to the computational complexity class NP. This class can be conveniently defined using the notion of a {\bf proof system}. A language (a set of binary strings) is in NP iff there is 
a {\bf p-bounded} proof system: a p-time decidable relation $P(x,y)$ (provability relation) 
and a 
constant $c \geq 1$ such that for any binary string $u$:
\begin{itemize}
	\item $u \in L$ iff $\exists w (|w| \le |u|^c)\ 
	P(u, w)$.
	
\end{itemize}
The existence of such a p-bounded system for the solvability problem is easy to see: $w$ is an assignment to variables 
and the provability relation $P$ simply says that all polynomials in $\fa$ vanish under the assignment. In fact, it is well-known that the problem is equally hard as the Boolean satisfiability problem and hence it is NP-complete.

\bigskip

Our research in this paper 
is motivated by the fundamental problem of proof complexity
whether NP is closed under complementation (the NP vs. coNP problem). 
By the NP-completeness mentioned above, this problem is thus equivalent to
the question whether one can also find a p-bounded proof system for the unsolvability of 
(\ref{8.1.23a}). More precisely, is there a p-time decidable relation $P(x,y)$ and a 
constant $c \geq 1$ such that for any $\cal F$:
\begin{itemize}
	\item $\cal F$ is unsolvable iff $\exists w (|w| \le len({\cal F})^c)\ P({\cal F}, w)$
	
\end{itemize}
where $len({\cal F})$ denotes the bit size of the representation of $\fa$.

An example of a proof system (but not p-bounded) for the unsolvability
is based on Hilbert's  Nullstellensatz and was introduced 
into proof complexity in \cite{BIKPP}. 
An {\bf NS-refutation} of (the solvability of) a set of polynomials $\cal F$ with variables $Var({\cal F})$ and
containing all polynomials $x^2 - x$ for all $x \in Var({\cal F})$ is a tuple of polynomials
$h_f \in \fp[Var(\fa)]$, for $f \in {\cal F}$, such that
$$
\sum_{f \in {\cal F}} \ h_f \cdot f\ =\ 1
$$
holds in $\fp[Var({\cal F})]$. The {\bf degree of the refutation} is 
$$
\max_{f \in {\cal F}} deg(h_f f)\ .
$$
It turns out that in natural situations the degree is the main measure we should care for.
Namely, in order to show some asymptotic lower bound it often suffices 
to consider finite sets of polynomials
$\fn$ over $\fp$ containing the polynomials $x^2-x$ for $x \in Var(\fn)$
such that for some constant $c \geq 1$ and all $n \geq 1$:
\begin{equation} \label{e2}
	|\fn|\le n^c\ ,\ N:= |Var(\fn)| \le n^c\ \mbox{ and }\ 
	\max \{deg(f)\ |\ f \in \fn\} \le c\ \ .
\end{equation}
A degree $d$ NS-refutation of $\fn$ has thus bit size
$O(|\fn| N^d) = O(n^{2c+d})$ which is polynomial in the bit size $len(\fn) = O(n^{3c})$
of $\fn$ iff $d$ is bounded by a constant. In other words, to show that the NS proof system is not p-bounded one needs a non-constant lower bound for $d$ as $n$ grows\footnote{For some other questions in proof complexity stronger lower bounds are needed, e.g. of the form $(\log n)^{\omega(1)}$ or even $n^{\Omega(1)}$.}. 
Such lower bounds are actually known, cf. \cite{BIKPP,BIKPRS} and references therein.

\bigskip

In this paper we study {\bf extended} NS proof systems which are stronger than the NS proof system.
The general idea of these proof systems is the following. 
Before looking for an NS refutation of $\fa \subseteq \rx$
extend it by a finite set $\ea \subseteq \rr$ of polynomials (to be called extension polynomials) in possibly more variables (to be called extension variables) than those of $\fa$, containing polynomials $r^p - r$ for all new variables $r$, 
and having the following {\bf soundness property}:
\begin{itemize}
	\item Any $0,1$ assignment $\overline a$ to variables $\overline x$
	can be appended by an $\fp$-assignment $\overline b$ to variables $\overline r$
	such that for all $g \in \ea$ it holds that $g(\overline a, \overline b) = 0$.
\end{itemize}
Then find an NS-refutation of $\fa \cup \ea$. Such proof may allow for a smaller degree than a mere NS refutation of $\fa$. However, a subtle point is that to make it a proof system one has to be able to check the soundness property of  
$\ea$ in p-time. This is a non-trivial requirement and for the two extended NS systems we consider in Section \ref{systems} this is enforced by syntactic requirements on sets $\ea$.

\medskip

In Section \ref{eval} we define a notion of pseudo-solutions of $\fa$ and prove that the existence of pseudo-solutions with suitable parameters
implies lower bounds for two extended NS proof systems ENS and UENS defined in \cite{BIKPRS}. 
In Section \ref{finitary} we give a combinatorial example of $\fa$ and candidate pseudo-solutions based on the pigeonhole principle. In Section \ref{logic} we 
discuss a motivation for extended
NS proof systems coming originally from logic (proof complexity).
The reader can find more proof complexity background in \cite{prf,BIKPRS}.

\smallskip

\noindent
{\bf A convention}: all logarithms are base $2$.

\section{Proof systems ENS and UENS} \label{systems}
Both examples of extended NS proof systems in this section, to be denoted ENS and UENS, are from \cite{BIKPRS}.

\begin{definition} [\cite{BIKPRS}] \label{2.1}
	{\ }
	
	\begin{enumerate}
		
		\item Let $\overline g = g_1, \dots, g_m$ be polynomials in any variables over $\ff_p$ and let $h \geq 1$ be a parameter.
		For $i \le m$ define polynomials:
		\begin{equation} \label{e4}
			E_{i,\overline g}\ :=\ g_i \cdot \Pi_{u \le h} (1 - \sum_{j \le m} r_{u j} g_j)
		\end{equation}
		where $r_{u j}$ are new {\bf extension variables} common to all $i\le m$.
		The polynomials $E_{i,\overline g}$ are called the
		{\bf extension polynomials of accuracy $h$} corresponding to $\overline g$.
		
		\item A set $\cal E$ of extension polynomials 
		can be {\bf stratified into $\ell$ levels} iff $\cal E$ can be partitioned
		as ${\cal E}_1 \cup \dots \cup {\cal E}_\ell$ where:
		
		\begin{itemize}
			
			\item If $E_{i, \overline g} \in {\cal E}_t$, some $t$ and $i$, then also
			all companion polynomials $E_{j, \overline g}$ are in ${\cal E}_t$, all $j \le m$.
			
			\item Variables in the polynomials $g_j$ in $E_{i, \overline g}$
			in ${\cal E}_1$ are among $Var({\cal F}_n)$
			and no extension variable from $E_{i, \overline g}$
			occurs among $Var({\cal F}_n)$ or in other extension polynomials in ${\cal E}_1$
			except in the companion polynomials.
			
			\item Variables in the polynomials $g_j$
			in the axioms $E_{i,\overline g}$ in ${\cal E}_{t+1}$, $1 \le t < \ell$, 
			are among the variables occurring in $Var({\cal F}_n \cup \bigcup_{s \le t} {\cal E}_s)$ (including 
			the extension variables from these levels)
			and no extension variable from $E_{i,\overline g}$ does occur among them or in other extension polynomials 
			in ${\cal E}_{t+1}$ except in the companion polynomials.
			
		\end{itemize}
	\item For a set $\cal E$ of extension polynomials put
	\begin{itemize}
		\item ${\cal R}({\cal E})$ to be the
			set of polynomials $r^p - r$, for all extension variables $r$ occurring in $\cal E$.
	\end{itemize}
	\end{enumerate}
\end{definition}
An {\bf ENS-refutation} of $\fa$ consists of a triple $(h, \ea, L)$, where 
$h \geq 1$ is its accuracy, $\ea$ is a set satisfying 
Definition \ref{2.1} and $L$ is an NS-refutation of $\fa \cup \ea \cup \re$. Its degree
is the degree of $L$.

\begin{lemma} \label{25.1.23a}
If $\fa$ has an ENS-refutation of any level $\ell \geq 1$ and any accuracy $h \geq 1$ then it is unsolvable. 
\end{lemma}

\prf

Let $\overline a$ be an assignment to variables $\overline x$.
Take any extension polynomial 
$E_{i,\overline g}$ in the first of $\ea$ and its factor
$(1 - \sum_{j \le m} r_{1 j} g_j)$. If all $g_j$ vanish under $\overline a$ 
put $r_{1 j}:=0$. If some $g_j$ does not vanish, put $r_{1 j'}:=0$ for $j' \neq j$ and
$r_{1 j}:= g_j(\overline a)^{-1}$. Further put $r_{i j}:=0$ for all $i > 1$ and any $j$.

This evaluation will also kill all companion polynomials $E_{i,\overline g}$, all $i$.
After you treat analogously other first level extension polynomials, move to the second level, etc. Finally note that polynomials in $\re$ vanish for any $\overline b$.

\qed

Extension polynomials were introduced in \cite{BIKPRS} 
so that the proof system can simulate in a sense
Boolean combinations of equations. Namely, assuming
that the equations $E_{i, \overline g} = 0$ hold for all $i \le m$,
we can express the truth value of the disjunction 
$$
\bigvee_{i \le m} (g_i \neq 0)
$$ 
by a polynomial
$$
{\sf disj}_{\overline g} \ :=\ 1\ -\ \Pi_{u \le h} (1 - \sum_{j \le m} r_{u j} g_j)
$$
of degree $h (1 + \max_i(deg(g_i))$ which may be much smaller, 
if $h$ is small, than $(p-1)\sum_j deg(g_i)$ required by the obvious definition $1 - \Pi_i (1 - g^{p-1}_i)$. 
Clearly if all $g_j = 0$ also ${\sf disj}_{\overline g} = 0$ 
and if $g_i \neq 0$ for some $i$ then, by $E_{i, \overline g} = 0$,
${\sf disj}_{\overline g} = 1$.
We shall discuss the motivation for introducing ENS in more detail in Section \ref{logic}.

\medskip

Extension polynomials can be also introduced in an unstructured form. A motivation for this
construction in \cite{BIKPRS} came from Boolean complexity.

\begin{definition}[{\cite{BIKPRS}}] \label{14.1.23a}
	An unstructured extension polynomial of accuracy $h \geq 1$ 
	is any polynomial from $\rr$ of the form 
$$
(g_1 - r_1) \cdot \dots \cdot (g_h - r_h)
$$
where no $r_i$ occurs in any $g_j$ (but may occur in other such unstructured extension
polynomials in $\ea$). 
\end{definition}

\begin{lemma} \label{15.1.23c}
Let $h \geq 1$ and let $\ea$ be a set of unstructured extension polynomials of accuracy
	$h$ such that $|\ea| < e^{h/p}$. Assume that $\fa \cup \ea \cup \re$ has an NS-refutation. Then
	$\fa$ is unsolvable.
\end{lemma}

\prf

For a given assignment $\overline a$ to variables $\overline x$ choose an assignment 
$\overline b$ to $\overline r$ uniformly at random from all $\fp$-assignments. 
The probability that one extension polynomial 
$(g_1 - r_1) \cdot \dots \cdot (g_h - r_h)$ does not vanish under $\overline a, \overline b$
is at most $(1 - 1/p)^h < e^{-h/p}$. 
Using the hypothesis that $|\ea| < e^{h/p}$ it follows by averaging
that there is some $\overline b$ such that all polynomials in $\ea$ vanish under $\overline a, \overline b$.
Finally note that polynomials in $\re$ vanish for any $\overline b$.

\qed

We can use this lemma to define a sound proof system:
a {\bf UENS-refutation} of $\fa \subseteq \rx$ is 
a triple $(h, \ea, L)$, where $h \geq 1$ is its accuracy, $\ea$ is a set satisfying satisfying
Definition \ref{14.1.23a} and also the size condition $|\ea| < e^{h/p}$,
and $L$ is an NS-refutation of $\fa \cup \ea \cup \re$. Its degree
is the degree of $L$.

UENS is at least as strong as ENS in the sense of the following statement.

\begin{lemma} [{\cite[Thm.6.14(1)]{BIKPRS}}] \label{15.1.23b}
For every extension polynomial $E_{i,\overline g}$ of accuracy $h \geq 1$
from Definition \ref{2.1}
there exists an unstructured extension polynomial $E'_{i,\overline g}$ of the same accuracy,
such that $E_{i,\overline g}$ can be expressed in $\rr$ as a linear combination of polynomials from the set
$E'_{i,\overline g}$, $x^2 -x$ for all $x$ and $r^p-r$ for all $r$.
of degree
at most $(ph + 1)\cdot \max_j deg(g_j) + h$.
\end{lemma}

Noting that the maximum degree of $E_{i,\overline g}$ for $j \le m$ is
$h(1 + \max_j deg(g_j)) +  \max_j deg(g_j)$ we get the following.

\begin{corollary} \label{21.1.23a}
	Assume that there is an ENS-refutation $h, \ea, L$
	of $\fa$ of degree $d$. Then there is an UENS-refutation $h, \ea', L'$
	of $\fa$  of degree $O(p d)$.
\end{corollary}

\section{Pseudo-solutions} \label{eval}

For $d \geq 0$ let $\rd$ be the $\fp$-vector space of polynomials from $\rx$ of degree
at most $d$.
We consider finite sets of polynomials
$\fn$ over $\fp$ in variables forming the set $Var(\fn)$, containing  
for each $x \in Var(\fn)$ the polynomial $x^2 - x$, and 
such that for some constant $c \geq 1$ conditions (\ref{e2}) are satisfied for $n \geq 1$.
We shall denote variables in $Var(\fn)$ simply $\overline x$.

A {\bf degree $d$ $\fn$-design} is any map
$$
\omega\ :\ \rd \rightarrow \fp
$$
satisfying the following three conditions:
\begin{enumerate}

\item For all $a \in \fp$: $\omega(a) = a$.

\item For all $g,h \in  \rd$: 
$$
\omega(g) + \omega(h) = \omega(g + h)\ .
$$ 

\item For any two polynomials $f \in \fn$ and $g \in \rd$, if $deg(f g) \le d$ then:
$$
\omega(f g) = 0\ .
$$

\end{enumerate}
This notion was defined by P.~Pudl\' ak in an informal electronic Prague-San Diego seminar
that run in early 1990s and it was used first in \cite{BIKPP} to establish the non-existence 
of constant degree NS-proofs of one modular counting principle from another, cf. also \cite[Sec.16.1]{prf} for an exposition of this and subsequent results.

An {\bf $\rd$-tree} $T$ is a finite 
$p$-ary tree whose each inner node (i.e. a non-leaf) is labeled by a query $g =\ ?$ for some
polynomial $g\in \rd$
and the $p$ outgoing edges are labeled by $g = a$, for all $a \in \fp$. The leaves
are labeled by elements of some set $I \neq \emptyset$. 
The {\bf height of $T$} is the maximum number of edges on a path from the root to
a leaf. An $\rd$-tree of height $\le e$ will be abbreviated as {\bf $(d,e)$-tree}.

Let $T$ be an $\rd$-tree whose leaves are labeled by elements of a non-empty set $I$.
The tree $T$ and any map $\omega : \rd \rightarrow \fp$ define a path in $T$ as follows: 
start at the root and answer the queries
by $\omega$. The label of a leaf the path so defined reaches
will be denoted $T(\omega)$. 
In particular, $T$ defines a map from the set of all degree $d$ $\fn$-designs to $I$.

We shall consider systems $\fn$ that have no solution in $\fp$. This means that no degree $d$ $\fn$-design
$\omega$ can be a homomorphism and there must be {\bf conflict pairs} of polynomials $g, g'\in \rd$, such that
$deg(g g') \le d$ and
$$
\omega(g) \cdot \omega(g') \neq \omega(g g')\ .
$$
The following notion was inspired by the model-theoretic construction in \cite[Cht.22]{k2} 
(it corresponds to the notion of sample space there).

\begin{definition}[pseudo-solutions] \label{1.1}
{\ }

For any $d, e \geq 0$ and $1 \geq  \gamma \geq 0$,
a $(d,e,\gamma)$-{\bf solution} of $\fn$
is a non-empty finite set $\Omega$ of degree $d$ $\fn$-designs $\omega$
such that for any $(d,e)$-tree $T$:
$$
\prob_{\omega \in \Omega}[T(\omega) \mbox{ is {\bf not} a conflict pair for $\omega$ }] \geq \gamma\ .
$$
A {\bf pseudo-solution} is a collective name for  
$(d,e,\gamma)$-solutions.

\end{definition}

\begin{theorem} \label{15.1.23a}
Let $(h, \ea, L)$ be an ENS- or UENS-refutation of $\fa$ of degree $d$. Let $S:= |\ea|$ and assume $e^{h/p} \geq 2 S^2$.
Then there is no $(d, h + \log S, S^{-1})$ pseudo-solution of $\fa$.

In particular, if $h$ is minimal such that $e^{h/p} \geq S^2$ holds then 
there is no $(d, (2p+1)\log S, S^{-1})$ pseudo-solution of $\fa$. 
\end{theorem}
The proof is summarized at the end of the section after establishing some lemmas.

Utilizing Corollary \ref{21.1.23a} we could prove the theorem only for UENS. However, because
of the importance of ENS for the logic problem discussed in Section \ref{logic} we shall
formulate the proof directly for ENS and avoid the detour via UENS; the argument can be
easily modified for UENS.
The idea is to transform in a sense the averaging argument for 
Lemma \ref{15.1.23c}
to the setting where we have $\fp$-linear maps $\omega : \rd \rightarrow \fp$ instead of
evaluation of variables by elements of $\fp$.

For $E_{i, \overline g}$ as in (\ref{e4}),
any of the $h+1$ polynomials 
\begin{equation} \label{9.7.21a}
g_i\ \mbox{ and }\  1 - \sum_{j \le m} r_{u j} g_j\ ,\ \mbox{ for }\ u \le h
\end{equation}
will be called a {\bf factor polynomial of $E_{i, \overline g}$}.

\begin{lemma} \label{2.3}
{\ }

Assume that $\Omega$ is a finite non-empty set of $\fp$-linear maps
from $\rd$ to $\fp$. Let
$\overline g = g_1, \dots, g_m$ be polynomials in variables $\overline x$
and let $E_{i,\overline g}$, $i \le m$, be the companion
extension polynomials of accuracy $h \geq 1$ using extension variables $r_{u j}$ as in (\ref{e4}).
Assume that all $E_{i,\overline g}$ have degree $\le d$.

Then for any $\omega \in \Omega$ and any $i \le m$, the probability over uniform random choices $\overline b$
that
$$
\omega(g'(\overline x, \overline b)) \neq 0
$$
holds for all factor polynomials $g'(\overline x, \overline r)$ of $E_{i, \overline g}$
is at most $e^{-h/p}$.
\end{lemma}

\prf

Let $\omega$ and $i \le m$ be fixed and fix at random all $b_{u j}$ with $u \le h$ and $j \neq i$.
Assume $\omega(g_i) \neq 0$.

Think of the remaining $h$ values 
$b_{u i}$, $u \le h$, as being chosen from $\fp$ uniformly and independently at random.
Using the $\fp$-linearity of $\omega$ we see that the probability to satisfy the inequality 
$$
\omega(\sum_{j \le m, j \neq i} b_{u j}g_j) + b_{u i}\omega(g_i) \neq 1\ .
$$
for one $u \le h$ is $\frac{p-1}{p}$. 
Hence we fail to violate it for all $u \le h$ with the 
probability at most $(1 - \frac{1}{p})^h \le e^{-h/p}$.

\qed

\bigskip

Having Lemma \ref{2.3} we may now proceed as level by level, bottom up, 
to find a suitable assignment of values from $\fp$ to all extension variables in $\cal E$. 
In particular, take one block of $m$ companion extension polynomials in ${\cal E}_1$
and choose (by an averaging argument) a suitable assignment to its extensions variables 
such that some of them fails to contain a factor polynomial that vanishes
for at most a fraction of $m \cdot e^{-h/p}$
of all $\omega$, and substitute it in the whole of $\cal E$. 
Then do the same consecutively for all blocks of companion polynomials in the first level. We lose altogether
at most $|{\cal E}_1| \cdot e^{-h/p}$ elements $\omega$. After we have substituted values from $\fp$
for all extension variables occurring in the first level, all polynomials $g_i$
in the second level ${\cal E}_2$ have only variables $\overline x$ left. 
Hence we can repeat the process. This yields the following statement.

\begin{lemma} \label{2.4a}
{\ }

Assume that $\Omega$ is a finite non-empty set of $\fp$-linear maps
from $\rd$ to $\fp$. Let
$\cal E$ be a leveled set of extension polynomials of accuracy $h$ of the form 
as in Def.\ref{2.1}, all of degree at most $d$.

Then there is an evaluation $\overline r := \overline b$ of the extension variables in $\cal E$
by values from $\fp$ 
such that for all $\omega \in \Omega$ but a fraction of $|{\cal E}| e^{-h/p}$ it holds that
\begin{itemize}

\item each extension polynomial in ${\cal E}$ has a factor polynomial $g'(\overline x, \overline r)$
such that
$$
\omega(g'(\overline x, \overline b)) = 0\ .
$$

\end{itemize}
\end{lemma}

\bigskip

\noindent
Now we are ready for 
{\bf Proof of Theorem \ref{15.1.23a}}.
Assume that ${\cal F}_n$ has an ENS-refutation $(h,\ea,L)$ of degree $d$ and that
$2 S^2 \le e^{h/p}$ holds where $S:= |\ea|$.
Hence:
\begin{equation} \label{20.7.21a}
|{\cal E}| e^{-h/p} \le  1/(2S)
\end{equation}
Note that, in particular, if $h$ was minimal to satisfy (\ref{20.7.21a}) then  
\begin{equation} \label{29.8.19b}
h  = O(\log S)\ .
\end{equation}
Let $L$:
$$
1\ =\ \sum_{f \in {\cal F}_n} h_f f\ +\ \sum_{E \in {\cal E}} p_E E\ +\ 
\sum_{g \in {\cal R}({\cal E})} q_g \cdot g
$$
be some NS-refutation of ${\cal F}_n \cup {\cal E} \cup {\cal R}({\cal E})$ of degree $d$.

\medskip

Assume that $\Omega$ is a $(d, h + \log S, S^{-1})$-solution of ${\cal F}_n$.
Substitute for all extension variables in $\cal E$
an evaluation $\overline b$ by elements of $\fp$ with the property stated in Lemma \ref{2.4a}.
In particular, 
after the substitution let $\Omega_0$ be the part of $\Omega$ consisting of $\omega$ 
satisfying the property in that lemma. Note that by the choice of $h$ 
\begin{equation} \label{29.8.19c}
|\Omega_0|\ \geq (1 - 1/(2S)) \cdot |\Omega|\ . 
\end{equation}

All polynomials $g \in {\cal R}({\cal E})$ vanish after the substitution 
and because maps $\omega \in \Omega$ are designs
also 
$\omega(\sum_{f \in {\cal F}_n} h_f(\overline b) f) = 0$ holds.
Hence we have that
\begin{equation} \label{3.6.18a}
1\ =\ \omega(\sum_{E \in {\cal E}} p_E(\overline b) E(\overline b))
\end{equation}
where the polynomials $p_E(\overline b), E(\overline b)$ have variables left only from
$Var({\cal F}_n)$.

\bigskip

We shall construct in two steps a particular $(d,e)$-tree $T^*$:

\begin{enumerate}

\item Use binary search to the sum in the right-hand side of (\ref{3.6.18a})
to get a $(d,\log S)$-tree $T_1$ 
that finds, for any $\omega \in \Omega$ a term $p_E E$ in the sum
which is
of the form $E = E_{i, \overline g}$ and $\omega(p_E E) \neq 0$.

\item To get $T^*$ extend $T_1$ as follows. 
At leaves labeled by $p_E E$ with $E = E_{i, \overline g}$, ask for values of all 
factor polynomials $1 - \sum_{j \le m} r_{u j} g_j$ for $u \le h$.
By the choice of $\Omega_0$, if $\omega \in \Omega_0$, 
at least one  gets by $\omega$ value $0$.
Hence we can write $p_E E$ as a product $g'g''$ with $\omega(g' g'') \neq 0 \wedge \omega(g') = 0$.
Label the corresponding leaf by $(g',g'')$.

\end{enumerate}
Note that the height $e$ of $T^*$ is $h + \log S$ 
and that for all $\omega \in \Omega_0$ it holds that $T^*(\omega)$ is a conflict pair for $\omega$.

\medskip

By the definition of pseudo-solutions we have that $T^*(\omega)$ can be a conflict pair for 
at most $(1 - 1/S)$-part of $\Omega$ but $\Omega_0$ is bigger. 
This is a contradiction and Theorem \ref{15.1.23a} is proved.

\qed

The relevance of the following problem is described in Section \ref{logic}.

\begin{problem} \label{15.1.23d}
Construct families $\fn$ (containing all polynomials $x^2 - x$ and obeying
(\ref{e2}) for some constant $c \geq 1$) and for any constant
$r \geq 1$ and $n >> 1$ an $((\log n)^r, r \log n, n^{-r})$ pseudo-solution
of $\fn$.	
\end{problem}

\section{A combinatorial example} \label{finitary}

The propositional formula $\php_n$:
$$
\bigvee_{i} \bigwedge_j \neg p_{ij}\ \vee\  
\bigvee_{i_1\neq i_2, j} (p_{i_1 j} \wedge p_{i_2 j})\ \vee\  
\bigwedge_{i, j_1 \neq j_2} (p_{i j_1} \wedge p_{i j_2})
$$
with $i, i_1, i_2$ ranging over $[n+1] = \{1, \dots, n+1\}$ and $j, j_1, j_2$ over $[n]$
is the most famous tautology in proof complexity, introduced by Cook and Reckhow \cite{CooRec}.
If it would be falsified by some evaluation $p_{i j} := a_{i j} \in \bits$ then the set
of all pairs $(i,j)$ such that $a_{i j} = 1$ would be the graph of an injective function
from $[n+1]$ into $[n]$. Hence the fact that $\php_n$ is logically valid is
equivalent to the pigeonhole principle.

The formula $\php_n$ requires long proofs in a number of proof systems
and its advanced variants may be hard for all proof systems (cf. \cite[Chpts.29-30]{k2} or
\cite[Chpt.19]{prf}). It has a polynomial size
(measured in its size $O(n^{3})$) proof in the ordinary propositional calculus
using the DeMorgan language and based on a finite number of axiom schemes and schematic inference rules,
a Frege system in the established terminology of \cite{CooRec}; this was proved by Buss \cite{Bus-php}.
In fact, it has a simple polynomial size proof in a $\mbox{TC}^0$-Frege system, a proof system operating
with bounded depth formulas using also the threshold connectives, and Buss's argument shows that
Frege systems do p-simulate these systems when such a connective is defined by a suitable DeMorgan formula.
One of the crucial results in proof complexity is Ajtai's theorem \cite{Ajt88} that proofs of $\php_n$ in
$\mbox{AC}^0$-Frege systems, 
subsystems of Frege systems operating with DeMorgan formulas of a bounded
depth, require super-polynomial number of steps (this was later improved to an exponential
lower bound in \cite{KPW,PBI}).

The negation of $\php_n$ can be reformulated (following \cite{BIKPP}) as the following 
unsolvable system $\pphp_n$
of polynomial equations over $\fp$ in variables $x_{i j}$, $i \in [n+1], j \in [n]$:
\begin{itemize}

\item $x_{i_1 j} \cdot x_{i_2 j} = 0$,  
for each $i_1 \neq i_2 \in [n+1]$ and $j \in [n]$.

\item $x_{i j_1} \cdot x_{i j_2} = 0$,  
for each $i \in [n+1]$ and $j_1 \neq j_2 \in [n]$.

\item $1 - \sum_{j \in [n]} x_{i j} = 0$, for each $i \in [n+1]$.

\item $x^2_{i j} - x_{i j}$, for all $i \in [n+1], j\in [n]$.

\end{itemize}
The left-hand sides of the equations in the first three items will be denoted $Q_{i_1, i_2; j}$,
$Q_{i; j_1, j_2}$ and $Q_i$, respectively.
We included the polynomials $x^2_{i j} - x_{i j}$ in $\pphp_n$ to conform with the requirement
put on systems ${\cal F}_n$ but it is easy to see that they are simple linear combinations of
the other polynomials:
$$
x^2_{i j} - x_{i j}\ =\ 
\sum_{k\in [n], k \neq j} Q_{i;j,k}\  -\ x_{i j} Q_i.
$$

NS-refutations of $\pphp_n$ need degree $n/2$ 
and that implies the following well-known fact. 

\begin{lemma} [\cite{Raz98,Buss98,BCEIP}]\label{10.7.21a}
For $n \geq 2$ there are degree $n/2$ $\pphp_n$-designs (over any $\fp$).
\end{lemma}

For a partial injective map $\rho : \subseteq [n+1] \rightarrow [n]$ define the {\bf restriction} of a polynomial
$g$ over $Var(\pphp_n)$ to be the polynomial to be denoted $g^\rho$ 
resulting from $g$ by the following partial substitution
of $0/1$ values for some variables:

\[ x_{i j}^{\rho} =  \left\{ \begin{array}{ll}
 1       &  \mbox{if $i \in dom(\rho) \wedge \rho(i)=j$} \\
 0       &  \mbox{if $i \in dom(\rho) \wedge \rho(i)\neq j$} \\
 0       &  \mbox{if $j \in rng(\rho) \wedge \rho^{(-1)}(j)\neq i$} \\
 x_{i j}  & \mbox{otherwise}
                    \end{array}
            \right. \]
Further define $D^{\rho} := [n+1] \setminus dom(\rho)$, 
$R^{\rho} := [n] \setminus rng(\rho)$ and
$n_{\rho} := |R^{\rho}| (= n - |\rho|)$.
Note that $x_{ij}^\rho = x_{ij}$ iff $(i,j) \in D^\rho \times R^\rho$.

If we apply a restriction $\rho$ to all polynomials in $\pphp_n$ we get a set 
$(\pphp_n)^\rho$ of
polynomials that expresses $\neg \php$ over $D^\rho$ and $R^\rho$ (plus the zero polynomial)
and that is isomorphic to $\pphp_{n_\rho}$. In particular, by Lemma \ref{10.7.21a} there are
degree $n_\rho/2$ designs for $(\pphp_n)^\rho$.

\begin{definition}
For $n \geq 1$ and $\epsilon > 0$ such that $n^\epsilon \geq 1$ let $\one$ be 
the set of all $\pphp_n$-designs $\omega$ that are defined as follows:
\begin{enumerate}

\item Pick 
\begin{enumerate}
\item a restriction $\rho : \subseteq [n+1] \rightarrow [n]$ with $n_\rho = n^\epsilon$,

\item a degree $n_\rho/2$-design $L$ for $(\pphp_n)^\rho$.

\end{enumerate}

\item For $g \in \phpd$ put:
$$
\omega(g)\ :=\ L(g^\rho)\ .
$$
\end{enumerate}
We shall write $\omega = (\rho_\omega, L_\omega)$ with $\rho_\omega$ and $L_\omega$ denoting the respective
$\rho$ and $L$ defining $\omega$.

\end{definition}

The next lemma interprets the task to find conflict pairs over $\one$ more combinatorially.

\begin{lemma}
Let $T$ be a $(d,e)$-tree over $\phpd$. Then there is $(d,e')$-tree $T'$ with 
$e'\le e + O(d \log n)$ such that for any $\omega = (\rho_\omega, L_\omega) \in \one$ if $T(\omega)$ is a conflict pair
for $\omega$ then $T'(\omega)$ is a pair $(i,j) \in D^\rho \times R^\rho$.
\end{lemma}

\prf

Assume $T$ found a conflict pair $g, g'$ for $\omega$. Write $g$ as an $\fp$-linear combination of monomials 
($\le n^{O(d)}$ of them) and use binary search to find one monomial $u$ such that
$$
\omega(u\cdot g') \neq \omega(u)\cdot \omega(g')\ .
$$
If $u$ is a product of variables $x_{i_1 j_1} \cdot \dots \cdot x_{i_t j_t}$
then by at most $2 t \le 2 d$ questions about values of 
$$
x_{i_s j_s}\ \mbox{ and }\ 
x_{i_1 j_1} \cdot \dots \cdot x_{i_s j_s} \cdot g'
$$
for $s = 1, \dots, t$ we find one variable $x_{i_s j_s}$ such that
$$
\omega(x_{i_s j_s})\cdot
\omega(x_{i_1 j_1} \cdot \dots \cdot x_{i_{s-1} j_{s-1}} \cdot g') \neq
\omega(x_{i_1 j_1} \cdot \dots \cdot x_{i_s j_s} \cdot g')\ .
$$
By the definition of restrictions this means that $x_{i_s j_s}^\rho = x_{i_s j_s}$. That is
$$
(i_s, j_s) \in D^\rho \times R^\rho\ .
$$

\qed

\begin{problem} \label{25.1.23b}
Do families $\pphp_n$ and sets $\one$ for some fixed $\epsilon > 0$
solve Problem \ref{15.1.23d}? 

Are $\one$, in fact,  
$(n^{\Omega(1)},n^{\Omega(1)}, 2^{-n^{\Omega(1)}})$ pseudo-solutions of $\pphp_n$, ?

\end{problem}

\section{Proof complexity motivation for ENS} \label{logic}

The problem to extend the lower bound for $\php_n$, or for any other formula for that matter,
from $\mbox{AC}^0$-Frege systems 
to $\mbox{AC}^0[p]$-Frege systems operating with formulas of a bounded 
depth in the DeMorgan language augmented by connectives counting modulo $p$, $p$ a prime, 
received a considerable attention over the last three decades.  
Proof-theoretically the most elegant definition of $\mbox{AC}^0$ proof systems
is using the formalism of sequent calculus LK (cf. \cite{kniha,prf})
but we shall stick to Frege systems: we will refer to 
a result from \cite{BIKPRS} and that used Frege systems. It is well known that the two formalisms yield
equivalent subsystems, except possibly for a change in the depth by an additive 
constant (cf. \cite[Chpt.3]{prf}).

Let $F$ be any Frege system in the DeMorgan language $0, 1, \neg, \vee,\wedge$. 
We shall denote by
$F(\mbox{MOD}_p)$ the proof system whose language extends the DeMorgan one 
by unbounded arity connectives $\mbox{MOD}_{p,i}$ for $p$ a prime and $i = 0, \dots, p-1$.
The formula $\mbox{MOD}_{p,i}(y_1, \dots, y_k)$ is true iff $\sum_j y_j \equiv i\ (\mbox{mod } p)$. 
The proof system has all Frege rules of $F$ accepted for all formulas of the extended language and, in
addition, the following set of $\mbox{MOD}_p$-axioms (cf. \cite[Sec.12.6]{kniha}):
\begin{itemize}
	\item $\mbox{MOD}_{p,0}(\emptyset)$
	
	\item $
	\neg \mbox{MOD}_{p,i}(\emptyset)\ ,\ \mbox{for}\ i=1,\dots, p-1
	$
	
	\item $
	\mbox{MOD}_{p,i}(\Gamma,\phi) \ \equiv\ [(\mbox{MOD}_{p,i}(\Gamma) \wedge
	\neg \phi) \vee (\mbox{MOD}_{p,i-1}(\Gamma) \wedge \phi)]
	$
	
	for $i = 0, \dots, p-1$, where $0-1$ means $p-1$ 
	and where $\Gamma$ stands for any sequence (possibly empty)
	of formulas.
\end{itemize}
The {\bf depth} of a constant or of an atom is $0$, the use of the negation or of any
of $\mbox{MOD}_{p,i}$ increases the depth
by $1$, and a formula formed from formulas $A_i$ none of which starts with $\vee$ (resp. with $\wedge$)
by a repeated use of $\vee$ (resp. of $\wedge$) has the depth $1$ plus the maximum depth
of formulas $A_i$. In particular, the depth of $\php_n$ is $3$. 
The subsystem of $F(\mbox{MOD}_p)$ allowed
to use only formulas of depth at most $\ell$ will be denoted by $F_\ell(\mbox{MOD}_p)$.

Any polynomial equation $f=0$ over $\ff_p$ is, in particular, also a
depth $2$ formula in the language of $F(\mbox{MOD}_p)$: monomials in $f$ can be defined by conjunctions
and $f=0$ is expressed using $\mbox{MOD}_{p,0}$ applied to them.
Therefore we can talk about $F(\mbox{MOD}_p)$-refutations of ${\cal F}_n$.
The following lemma is straightforward (cf. also \cite{BIKPP}).

\begin{lemma} \label{easy}
	The two formulations 
	$\neg \php_n$ and $\pphp_n$ can be derived form each other in $F_4(\mbox{MOD}_p)$
	by proofs of size $n^{O(1)}$.
\end{lemma} 
Having an $F_\ell(\mbox{MOD}_p)$-proof we can translate it into a proof using only low degree polynomials;
first define big conjunctions $\bigwedge$ via big disjunctions $\bigvee$ and $\neg$, and then 
translate bottom up all $\bigvee$ by ${\sf disj}$ as described after Lemma \ref{25.1.23a}, while systematically introducing all needed
extension polynomials. This yields the following result.

\begin{theorem} [{\cite[Thm.6.7(1)]{BIKPRS}}] \label{2.2}
	{\ }
	
	Let $\ell \geq 2$ be a constant.
	Let ${\cal F}_n$ be a set of polynomials obeying (\ref{e2}) and 
	containing all polynomials $x^2 - x$ for all $x \in Var({\cal F}_n)$, and 
	assume it has an $F_\ell(\mbox{MOD}_p)$-refutation with $k$ steps. 
	Let $h \geq 1$ be any parameter.
	
	Then for any $h \geq 1$ 
	there exists a set $\cal E$ of $S:=k^{O(1)}$ extension polynomials of accuracy $h$
	stratified into $\ell + O(1)$ levels such that ${\cal F}_n \cup {\cal E} \cup {\cal R}({\cal E})$ 
	has an $NS$-refutation of degree at most 
	\begin{equation} \label{e5}
		(O(1) + \log k)(h+1)^{O(1)}\ .
	\end{equation}
\end{theorem}
(The constant $O(1)$ in term $k^{O(1)}$ and the additive constant $O(1)$ depend on the underlying Frege system $F$ only but the $O(1)$ in the exponent in (\ref{e5}) 
depends also on $\ell$: the construction underlying the proof in \cite{BIKPRS} yields 
the bound $O(\ell)$.)

\bigskip

Should there be refutations of $\fn$ in any $F_\ell(\mbox{MOD}_p)$, 
fixed $\ell \geq 2$, with $k = n^{O(1)}$ steps, we can choose minimal $h \geq 1$
such that $e^{h/p} \geq 2 S^2$, $S$ provided by Theorem \ref{2.2}. Such $h$ is $O(\log n)$
and hence we would get ENS-refutations $(h, \ea, L)$ of degree 
$$
d \le (O(1) + \log k)(h+1)^{O(1)} \le (\log n)^c
$$
where $c$ depends on $\ell$.
This connection\footnote{An earlier reduction in \cite{Kra-fdp} reduced the lengths-of-proofs problem for $\mbox{AC}^0[p]$-Frege systems to a different and seemingly harder problem about algebraic search trees.}
to ENS, and Theorem \ref{15.1.23a} , is the reason why we are interested in Problem \ref{25.1.23a}. 

Note that should it be answered by $\pphp_n$ (Problem \ref{25.1.23b})
it would imply also a  lower bound for UENS-refutations of $\pphp_n$. 
This would disprove a causal remark\footnote{The construction in \cite{BIKPRS} seems to yield
a non-trivial degree upper bound on ENS proof when transforming Frege proofs of 
logical depth $o(\log n)$ only.} 
at the end of \cite{BIKPRS} that size $S$ Extended Frege (EF) refutations can be transformed into $(\log S)^{O(1)}$ degree UENS-refutations as it is well-known that
EF admits p-size proofs of $\php$, cf. \cite{CooRec}.

Note that it is also unknown if EF simulates UENS.

\bigskip

\noindent
{\large {\bf Acknowledgments.}} I am indebted to the anonymous referee for
a number of suggestions.

\end{document}